%% file: wmain.tex
\documentclass[12pt,leqno,twoside]{article} 
\usepackage{bezier,amsmath,amssymb}
\input{rogg_kue.sty}
\newcommand{\mufat}{\boldsymbol{\mu}}
\newcommand{\lafat}{\boldsymbol{\lambda}}

\begin{document}
\title{On coefficient valuations of Eisenstein polynomials}
\author{M.\ K\"unzer, E.\ Wirsing}
\maketitle

\begin{small}
\begin{quote}
\begin{center}{\bf Abstract}\end{center}\vspace*{2mm}
Let $p\geq 3$ be a prime, let $n > m\geq 1$. Let $\pi_n$ be the norm of $\zeta_{p^n} - 1$ under
$C_{p-1}$, so that $\Z_{(p)}[\pi_n]|\Z_{(p)}$ is a purely ramified extension of discrete valuation rings of degree $p^{n-1}$. The minimal polynomial of $\pi_n$ 
over $\Q(\pi_m)$ is an Eisenstein polynomial; we give lower bounds for its coefficient valuations at $\pi_m$. The function field analogue, as introduced by 
Carlitz and Hayes, is studied as well.
\end{quote}
\end{small}

\renewcommand{\thefootnote}{\fnsymbol{footnote}}
\renewcommand{\thefootnote}{\arabic{footnote}}

\begin{footnotesize}
\renewcommand{\baselinestretch}{0.7}
\parskip0.0ex
\tableofcontents
\parskip1.2ex
\renewcommand{\baselinestretch}{1.0}
\end{footnotesize}

\input{wsec0}         
\input{wsec0_5}       
\input{wsec1}         
\input{wsec2}         
\input{wsec3}         
\input{wsec4}         
\input{wsec5}         
\input{wref}          
\end{document}

%% file: wsec0.tex
\setcounter{section}{-1}

\section{Introduction}

\subsection{The problem}
\label{SubsecSit}
Suppose given a purely ramified extension of discrete valuation rings $S|R$, with maximal ideals generated by $s\in S$ and $r\in R$, respectively. In particular, 
$S = R[s]$. Let $L|K$ be the corresponding extension of the fields of fractions, and write $l = [L:K]$. The minimal polynomial $\mu_{s,K}(X)\in R[X]$ of $s$ over 
$K$ is an Eisenstein polynomial, that is, the valuation at $r$ of its non-leading coefficients is $\geq 1$, and the valuation of its constant term equals $1$. 

Thus $S/rS\iso (R/rR)[X]/(X^l)$ is completely known. To obtain information about 
\[
S/r^m S \;\iso\; (R/r^m R)[X]/(\mu_{s,K}(X))
\]
for $m\geq 1$, however, we need to know better estimates for the coefficient valuations of the minimal polynomial $\mu_{s,K}(X)$. In short: how eisensteinian is 
it really?

The objective of this note is to give lower valuation bounds for the coefficients of $\mu_{s,K}(X)$ for certain subextensions $S|R$ of cyclotomic extensions, 
both in the number field and in the function field case. Moreover, the method we use enables us to relate the minimal polynomials for $T|R$ and for $S|R$
for iterated extensions $R\tm S\tm T$ of discrete valuation rings therein.

As an application, we mention the Wedderburn embedding of the twisted group ring (with trivial $2$-cocycle)
\[
S\wr C_l\;\hraa{\omega}\; \End_R S \;\iso\; R^{l\ti l}\; ,
\]
where we assume $L|K$ to be galois with cyclic Galois group $C_l$. We have
\[
\omega(S\wr C_l) \;\tm\; \Lambda \; :=\; \{ f\in \End_R S \, :\, \mb{$f(s^i S) \tm s^i S$ for $i\geq 0$}\} \; \tm\; \End_R S\; . 
\] 
The image of $s$ is the companion matrix of $\mu_{s,K}(X)$. To give a description of $\omega(S\wr C_l)$, it is convenient to be able to replace the image 
of $s$ by the companion matrix of the `optimal' Eisenstein polynomial $X^l - r$ modulo $s^j\Lambda$ for a suitable $j > 1$. 

In this article, however, we restrict our attention to the minimal polynomial itself.

\subsection{Results}

\subsubsection{The number field case}

Let $p\geq 3$ be a prime, and let $\zeta_{p^n}$ denote a primitive $p^n$th root of unity over $\Q$ in such a way that $\zeta_{p^{n+1}}^p = \zeta_{p^n}$ for all 
$n\geq 1$. Put $F_n = \Q(\zeta_{p^n})$ and let $E_n = \mb{Fix}_{C_{p-1}} F_n$, so $[E_n:\Q] = p^n$. Letting
$\;
\pi_n \; =\; \Nrm_{F_n|E_n}(\zeta_{p_n} - 1) \; =\; \prod_{j\in [1,p-1]} (\zeta_{p^n}^{j^{p^{n-1}}} - 1)\; ,
$
we have $E_n = \Q(\pi_n)$. In particular, $E_{m+i} = E_m(\pi_{m+i})$ for $m,i\geq 1$. We fix $m$ and write
\[
\mu_{\pi_{m+i},\,E_m}(X)\; =\; \sum_{j\in [0,p^i]} a_{i,j} X^j \; =\; X^{p^i} + \Big(\sum_{j\in [1,p^i - 1]} a_{i,j} X^j\Big) - \pi_m \;\in\; \Z[\pi_m][X]\; .
\]

\pagebreak[3]

{\bf Theorem} (\ref{Th11}, \ref{CorDiff}, \ref{CorEis13}).
{\it
\begin{itemize}
\item[{\rm (i)}] We have $p^i\; |\; j a_{i,j}$ for $j\in [0,p^i]$. 
\item[{\rm (i$'$)}] If $j < p^i (p-2)/(p-1)$, then $p^i\pi_m\; |\; j a_{i,j}$.
\item[{\rm (ii)}] We have $a_{i,j} \con_{p^{i+1}} a_{i + \be,p^\be j}$ for $j\in [0,p^i]$ and $\be\geq 1$. 
\item[{\rm (ii$'$)}] If $j < p^i (p-2)/(p-1)$, then $a_{i,j} \con_{p^{i+1}\pi_m} a_{i+\be,p^\be j}$ for $\be\geq 1$.
\item[{\rm (iii)}] We have a unit $a_{i,p^i - (p^i - p^\be)/(p-1)} \,/\, p^{i-\be} \;\in\; \Z_{(p)}[\pi_m]^\ast$ for $\be\in [0,i-1]$.
\item[{\rm (iv)}] We have
$\;
\mu_{\pi_n,\sQ}(X)\; \con_{p^2} X^{p^{n-1}} + pX^{(p-1)p^{n-2}} - p
$
for $n\geq 2$.
\end{itemize}
}

Assertion (iv) requires a computation of a trace. Such trace computations can be reformulated in terms of sums of $(p-1)$th roots of unity in $\Q_p$ 
(\ref{PropFT2}). Essentially, one has to count the number of subsets of $\mufat_{p-1}\tm\Q_p$ of a given cardinality whose sum is of a given valuation at $p$. 
Not being able to go much beyond this reformulation, this seems to be a problem in its own right -- see e.g.\ (\ref{ExEis7_5}).

To prove (i, i$'$, ii, ii$'$), we proceed by induction. Assertions (i, i$'$) also result from the different 
$\Dfk_{\sZ_{(p)}[\pi_{m+i}]|\sZ_{(p)}[\pi_m]} = \left(\mu'_{\pi_{m+i},E_m}(\pi_{m+i})\right) = \left(p^i \pi_{m+i}^{p^i - 1 - (p^i - 1)/(p - 1)}\right)$. 
Moreover, using (ii), this different argument yields (iii). In the function field case below, this argument for (i, i$'$) fails, however, and we have to 
resort to induction.

Suppose $m = 1$. Let us call an index $j\in [1,p^i - 1]$ {\it exact,} if either $j < p^i (p-2)/(p-1)$ and $p^i\pi_m$ exactly divides $j a_{i,j}$, or 
$j \geq p^i (p-2)/(p-1)$ and $p^i$ exactly divides $j a_{i,j}$. If $i = 1$ and e.g.\ $p\in \{3, 19, 29, 41\}$, then all indices $j\in [1,p - 1]$ are exact. 
If $i\geq 2$, we might ask whether the number of non-exact indices $j$ asymptotically equals $p^{i-1}$ as $p\to\infty$. 

\subsubsection{The function field case}

Let $p\geq 3$ be a prime, $\rho\geq 1$ and $r = p^\rho$. We write $\Zl = \Fu{r}[Y]$ and $\Ql = \Fu{r}(Y)$. We want to study
a function field analogue over $\Ql$ of the number field extension $\Q(\zeta_{p^n})|\Q$. Since $1$ is the only $p^n$th root of unity in an algebraic closure 
$\b\Ql$, we have to proceed differently, following Carlitz {\bf\cite{Ca38}} and Hayes {\bf\cite{Ha74}}. First of all, the power operation of $p^n$ on 
$\b\Q$ becomes replaced by a module operation of $f^n$ on $\b\Ql$, where $f\in\Zl$ is an irreducible polynomial. The group of $p^n$th roots of unity
\[
\mufat_{p^n}\=\{\xi\in\b\Q : \xi^{p^n} = 1\}\;\, 
\]
becomes replaced by the annihilator submodule 
\[
\lafat_{f^n}\=\{\xi\in\b\Ql : \xi^{f^n} = 0\}\; .
\]
Instead of choosing a primitive $p^n$th root of unity $\zeta_{p^n}$, i.e.\ a $\Z$-linear generator of that abelian group, we choose a $\Zl$-linear generator 
$\theta_n$ of this $\Zl$-submodule. A bit more precisely speaking, the element $\theta_n\in\b\Ql$ plays the role of $\th_n := \zeta_{p^n} - 1\in\b\Q$. Now 
$\Ql(\theta_n)|\Ql$ is the function field analogue of $\Q(\th_n)|\Q$. See also {\bf\cite[\rm sec.\ 2]{Go83}}.

To state the result, let $f(Y)\in \Zl$ be a monic irreducible polynomial and write $q = r^{\deg f}$. 
Let $\xi^Y := Y\xi + \xi^r$ define the $\Zl$-linear Carlitz module structure on an algebraic closure $\b\Ql$, and choose a $\Zl$-linear generator $\theta_n$ of 
$\mb{ann}_{f^n}\b\Ql$ in such a way that $\theta_{n+1}^f = \theta_n$ for all $n\geq 1$. We write $\Fl_n = \Ql(\theta_n)$ and have 
$\Gal(\Fl_n/\Q) \iso (\Zl/f^n)^\ast$. Letting $\El_n = \mb{Fix}_{C_{q-1}} \Fl_n$, we get $[\El_n:\Ql] = q^n$. Denoting
$\;
\varpi_n \; =\; \Nrm_{\Fl_n|\El_n}(\theta_n) \; =\; \prod_{e\in (\Zl/f)^\ast} \theta_n^{e^{q^{n-1}}}\; ,
$
we have $\El_n = \Ql(\varpi_n)$. In particular, $\El_{m+i} = \El_m(\varpi_{m+i})$ for $m,i\geq 1$. We fix $m$ and write
\[
\mu_{\varpi_{m+i},\,\El_m}(X)\; =\; \sum_{j\in [0,q^i]} a_{i,j} X^j \; =\; X^{q^i} + \Big(\sum_{j\in [1,q^i - 1]} a_{i,j} X^j\Big) - \varpi_m 
\;\in\; \Zl[\varpi_m][X]\; .
\]
Let $v_q(j) := \min\{ \alpha\geq 0 : \; q^\alpha\, |\, j \; \}$.

{\bf Theorem} (\ref{Th11a}, \ref{CorDiff2}, \ref{CorBU}). 
{\it
\begin{itemize}
\item[{\rm (i)}] We have $f^{i - v_q(j)}\; |\; a_{i,j}$ for $j\in [0,q^i]$. 
\item[{\rm (i$'$)}] If $j < q^i (q-2)/(q-1)$, then $f^{i - v_q(j)} \varpi_m\; |\; a_{i,j}$.
\item[{\rm (ii)}] We have $a_{i,j} \con_{f^{i+1}} a_{i+\be,q^\be j}$ for $j\in [0,q^i]$ and $\be\geq 1$. 
\item[{\rm (ii$'$)}] If $j < q^i (q-2)/(q-1)$, then $a_{i,j} \con_{f^{i+1}\varpi_m} a_{i+\be,q^\be j}$ for $\be\geq 1$.
\item[{\rm (iii)}] We have a unit $a_{i,q^i - (q^i - q^\be)/(q-1)} \,/\, f^{i-\be} \;\in\; \Zl_{(f)}[\varpi_m]^\ast$ for $\be\in [0,i-1]$.
\item[{\rm (iv)}] If $f = Y$, then
$\;
\mu_{\varpi_{m+i},\,\El_m}(X)\; \con_{Y^2} X^{q^i} + YX^{(q-1)q^{i-1}} - \varpi_m\; .
$
\end{itemize}
}

A comparison of the assertions (iv) in the number field case and in the function field case indicates possible generalizations; we do not know 
what happens for $\mu_{\pi_{m+i},E_m}(X)$ for $m\geq 2$ in the number field case; moreover, we do not know what happens for $f\neq Y$ in
the function field case.

\subsection{Notations and conventions}

\begin{footnotesize}
\begin{itemize}
\item[(o)] Within a chapter, the lemmata, propositions etc. are numbered consecutively.
\item[(i)] For $a,b\in\Z$, we denote by $[a,b] := \{c\in\Z\; :\; a\leq c\leq b\}$ the interval in $\Z$.
\item[(ii)] For $m\in\Z\ohne\{ 0\}$ and a prime $p$, we denote by $m[p] := p^{v_p(m)}$ the $p$-part of $m$, where $v_p$ denotes the valuation of an integer at $p$.
\item[(iii)] If $R$ is a discrete valuation ring with maximal ideal generated by $r$, we write $v_r(x)$ for the valuation of $x\in R\ohne\{ 0\}$ at $r$, i.e.\ 
$x/r^{v_r(x)}$ is a unit in $R$. In addition, $v_r(0) := +\infty$.
\item[(iv)] Given an element $x$ algebraic over a field $K$, we denote by $\mu_{x,K}(X)\in K[X]$ the minimal polynomial of $x$ over $K$.
\item[(v)] Given a commutative ring $A$ and an element $a\in A$, we sometimes denote the quotient by $A/a := A/aA$ --- mainly if $A$ plays the role of a base ring. For 
$b,c\in A$, we write $b\con_a c$ if $b - c\in aA$. 
\item[(vi)] For an assertion $X$, which might be true or not, we let $\{ X\}$ equal $1$ if $X$ is true, and equal $0$ if $X$ is false.
\end{itemize}
\end{footnotesize}

%% file: wsec0_5.tex
\fbox{Throughout, let $p\geq 3$ be a prime.}

\section{A polynomial lemma}

We consider the polynomial ring $\Z[X,Y]$.

\begin{Lemma}
\label{Lem2}
For $k\geq 1$, we have 
$\;
(X + pY)^k\;\con_{k[p]\cdot p^2 Y^2}\; X^k + k X^{k-1} pY \; .
$

\rm
Since $\smatze{k}{j} =  k/j \cdot \smatze{k-1}{j-1}$, we obtain for $j\geq 2$ that
\[
\begin{array}{rcl}
v_p(p^{\, j}\smatze{k}{j}) 
& \geq & j + v_p(k) - v_p(j) \\
& \geq & v_p(k) + 2\; , \\ 
\end{array}
\]
where the second inequality follows from $j\geq 2$ if $v_p(j) = 0$, and from $j\geq p^{v_p(j)}\geq 3^{v_p(j)}\geq v_p(j) + 2$ if $v_p(j)\geq 1$.
\end{Lemma}

\begin{Corollary}
\label{Cor2a}
For $k\geq 1$, we have 
$\;
(X + pY)^k \;\con_{k[p]\cdot p Y} X^k\; . 
$
\end{Corollary}

\begin{Corollary}
\label{Cor2aa}
For $x,y\in\Z$ and $l\geq 1$ such that $x\con_{p^l} y$, and for $k\geq 1$, we have
$\;
x^k\;\con_{k[p]\cdot p^l}\; y^k\; .
$
\end{Corollary}

\begin{Corollary}
\label{Cor2b}
For all $\alpha,\beta\geq 0$, we have
$\;
(X+Y)^{p^{\beta +\alpha}}\;\con_{p^{\alpha + 1}} (X^{p^\beta } + Y^{p^\beta})^{p^\alpha}\; .
$

\rm
This being true for $\alpha = 0$, the assertion follows since $f(X,Y)\con_p g(X,Y)$ implies that $f(X,Y)^{p^\alpha}\con_{p^{\alpha + 1}} g(X,Y)^{p^\alpha}$ by 
(\ref{Cor2a}), where $f(X,Y),\, g(X,Y)\in\Z[X,Y]$.
\end{Corollary}

%% file: wsec1.tex
\section{Consecutive purely ramified extensions}

\subsection{Setup}

Let $T|S$ and $S|R$ be finite and purely ramified extensions of discrete valuation rings, of residue characteristic $\ch R/rR = p$. The maximal ideals of $R$, $S$ and $T$ are 
generated by $r\in R$, $s\in S$ and $t\in T$, and the fields of fractions are denoted by $K = \fracfield R$, $L = \fracfield S$ and $M = \fracfield T$, respectively.  
Denote $m = [M:L]$ and $l = [L:K]$. We may and will assume $s = (-1)^{m+1}\Nrm_{M|L}(t)$ and $r = (-1)^{l+1}\Nrm_{L|K}(s)$. 

We have $S = R[s]$ with
\[
\mu_{s,K}(X) \; =\; X^l + \Big(\sum_{j\in [1,l - 1]} a_j X^j\Big) - r \;\in\; R[X]\; ,
\]
and $T = R[t]$ with
\[
\mu_{t,K}(X) \; =\; X^{lm} + \Big(\sum_{j\in [1,lm - 1]} b_j X^j\Big) - r \;\in\; R[X]\; .
\]
Cf.\ {\bf \cite[\rm I.\S 7, prop.\ 18]{Se62}}. The situation can be summarized in the diagram
\begin{center}
\begin{picture}(250,550)
\put(-105,   0){$rR\tm R$}
\put(  10,  50){\line(0,1){130}}
\put(-210, 202){$sS\tm S = R[s]$}
\put(  10, 250){\line(0,1){130}}
\put(-330, 402){$tT\tm T = S[t] = R[t]$}
\put(  50,  50){\line(1,1){40}}
\put(  50, 250){\line(1,1){40}}
\put(  50, 450){\line(1,1){40}}
\put( 100, 100){$K$}
\put( 110, 150){\line(0,1){130}}
\put( 120, 205){$\scm l$}
\put( 100, 300){$L$}
\put( 110, 350){\line(0,1){130}}
\put( 120, 405){$\scm m$}
\put(  95, 500){$M$}
\end{picture}
\end{center}

Note that $r\; |\; p$, and that for $z\in M$, we have $v_t(z) = m v_s(z) = ml v_r(z)$.

\subsection{Characteristic $0$}

In this section, we assume $\ch K = 0$. In particular, $\Z_{(p)}\tm R$.

\begin{Assumption}
\label{Ass3}\rm
Suppose given $x,y\in T$ and $k\in [1,l-1]$ such that
\begin{itemize}
\item[(i)] $p\; |\; y$ and $t^m \con_y s$,
\item[(ii)] $x\; |\; j a_j$ for all $j\in [1,l - 1]$, and
\item[(iii)] $xr\; |\; j a_j$ for all $j\in [1,k - 1]$.
\end{itemize}
\end{Assumption}

Put $c := \gcd(xys^{k-1},yls^{l-1})\in T$.

\begin{Lemma}
\label{Lem5}
Given {\rm (\ref{Ass3}),} we have $\; c\; |\; \mu_{s,K}(t^m)\;$.

\rm
We may decompose
\[
\begin{array}{rcl}
\mu_{s,K}(t^m)
& = & \mu_{s,K}(t^m) - \mu_{s,K}(s) \\
& = & (t^{ml} - s^l) + \Big(\sum_{j\in [1,k-1]} a_j (t^{mj} - s^j)\Big) + \Big(\sum_{j\in [k,l-1]} a_j (t^{mj} - s^j)\Big)\; . \\
\end{array}
\]
Now since $t^m = s + zy$ for some $z\in T$ by (\ref{Ass3}.i), we have 
\[
t^{mj} \;\auf{(\ref{Lem2})}{\con}_{j y^2} s^j + j s^{j-1} zy \;\con_{j s^{j-1} y} s^j
\]
for any $j\geq 1$, so that $s^{j-1}\; |\; r \; |\; p\; |\; y$ gives $t^{mj}\con_{js^{j-1}y} s^j$.

In particular, $y l s^{l-1}\; |\; t^{ml} - s^l$.

Moreover, $x y s^l\; |\;\sum_{j\in [1,k-1]} a_j (t^{mj} - s^j)$ by (\ref{Ass3}.iii).

Finally, $x y s^{k-1}\; |\; \sum_{j\in [k,l-1]} a_j (t^{mj} - s^j)$ by (\ref{Ass3}.ii). 
\end{Lemma}

The following proposition will serve as inductive step in (\ref{Prop7}).

\begin{Proposition}
\label{Prop6}
Given {\rm (\ref{Ass3}),} we have $t^{-j} c\; |\; b_j$ if $j\not\con_m 0$ and $t^{-j} c\; |\; (b_j - a_{j/m})$ if $j\con_m 0$ for $j\in [1,lm - 1]$.

\rm
From (\ref{Lem5}) we take 
\[
\sum_{j\in [1,lm - 1]} \left(b_j - \{ j\con_m 0\}\, a_{j/m}\right) t^j \; =\; -\mu_{s,K}(t^m) \;\con_c\; 0\; . 
\]
Since the summands have pairwise different valuations at $t$, we obtain
\[
\left(b_j - \{ j\con_m 0\}\, a_{j/m}\right) t^j \;\con_c\; 0
\]
for all $j\in [1,lm - 1]$.
\end{Proposition}

\subsection{As an illustration: cyclotomic polynomials}

\begin{quote}
\begin{footnotesize}
For $n\geq 1$, we choose primitive roots of unity $\zeta_{p^n}$ over $\Q$ in such a manner that $\zeta_{p^{n+1}}^p = \zeta_{p^n}$. We abbreviate 
$\th_n = \zeta_{p^n} - 1$.

We shall show by induction on $n$ that writing
\[
\mu_{\th_n,\sQ}(X) \; =\; \Phi_{p^n}(X+1)\; =\; \sum_{j\in [0,p^{n-1}(p-1)]} d_{n,j} X^j
\]
with $d_{n,j}\in\Z$, we have $p^{n-1}\; |\; j d_{n,j}$ for $j\in [0,p^{n-1}(p-1)]$, and even $p^n\; |\; j d_{n,j}$ for $j\in [0,p^{n-1}(p-2)]$.

This being true for $n = 1$ since $\Phi_p(X+1) = ((X+1)^p - 1)/X$, we assume it to be true for $n - 1$ and shall show it for $n$, where $n\geq 2$.
We apply the result of the previous section to $R = \Z_{(p)}$, $r = -p$, $S = \Z_{(p)}[\th_{n-1}]$, $s = \th_{n-1}$ and $T = \Z_{(p)}[\th_n]$, $t = \th_n$. 
In particular, we have $l = p^{n-2}(p-1)$ and $\mu_{s,K}(X) = \Phi_{p^{n-1}}(X+1)$; we have $m = p$ and $\mu_{t,L}(X) = (X + 1)^p - 1 - \th_{n-1}$; finally, 
we have $\mu_{t,K}(X) = \Phi_{p^n}(X+1)$.

We may choose $y = p\th_n$, $x = p^{n-2}$ and $k = p^{n-2}(p-2) + 1$ in (\ref{Ass3}). Hence $c = p^{n-1}\th_n^{p^n - 2 p^{n-1} + 1}$. By (\ref{Prop6}), we obtain 
that $p^{n-1}\th_n^{p^n - 2 p^{n-1} + 1 - j}$ divides $d_{n,j} - d_{n-1,j/p}$ if $j\con_p 0$ and that it divides $d_{n,j}$ if $j\not\con_p 0$. Since the 
coefficients in question are in $R$, we may draw the following conclusion. 
$$
\left\{
\mb{
\begin{tabular}{p{12cm}}
If $j\con_p 0$, then $p^n \; |\;  d_{n,j} - d_{n-1,j/p}$ if $j \leq p^{n-1}(p-2)$, \\
\hspace*{11.8mm} and $p^{n-1} \; |\;  d_{n,j} - d_{n-1,j/p}$ if $j > p^{n-1}(p-2)$; \\
if $j\not\con_p 0$, then $p^n \; |\;  d_{n,j}$ if $j \leq p^{n-1}(p-2)$, \\
\hspace*{11.5mm} and $p^{n-1} \; |\;  d_{n,j}$ if $j > p^{n-1}(p-2)$. \\
\end{tabular}
}
\right.
\leqno (\mb{I})
$$
By induction, this establishes the claim. 

Using (\ref{Cor2b}), assertion (I) also follows from the more precise relation
$$
\Phi_{p^n}(X+1) - \Phi_{p^{n-1}}(X^p + 1) \;\con_{p^n}\; X^{p^{n-1}(p - 2)} \left( (X^p + 1)^{p^{n-2}} - (X + 1)^{p^{n-1}}\right)
\leqno (\mb{II})
$$
for $n\geq 2$, which we shall show now. In fact,
\[
\begin{array}{cl}
   & \left(\ru{-1.5}(X+1)^{p^n} - 1\right)\left((X^p + 1)^{p^{n-2}} - 1\right) - \left((X^p+1)^{p^{n-1}} - 1\right)\left((X + 1)^{p^{n-1}} - 1\right) \\
\auf{\mb{\scr (\ref{Cor2b})}}{\con}_{\!\! p^n} & \left((X^p+1)^{p^{n-1}} - 1\right)\left((X^p + 1)^{p^{n-2}} - (X + 1)^{p^{n-1}}\right) \\
\auf{\mb{\scr (\ref{Cor2b})}}{\con}_{\!\! p^n} & X^{p^n}\left((X^p + 1)^{p^{n-2}} - (X + 1)^{p^{n-1}}\right) \\
\auf{\mb{\scr (\ref{Cor2b})}}{\con}_{\!\! p^n} & X^{p^{n-1}(p - 2)}\left((X^p + 1)^{p^{n-2}} - (X + 1)^{p^{n-1}}\right)
                                                 \left((X+1)^{p^{n-1}} - 1\right)\left((X^p+1)^{p^{n-2}} - 1\right) \\
\end{array}
\]
and the result follows by division by the monic polynomial 
\[
\left((X+1)^{p^{n-1}} - 1\right)\left((X^p+1)^{p^{n-2}} - 1\right)\; .
\]
Finally, we remark that writing $F_n(X) := \Phi_{p^n}(X+1) + X^{p^n - 2p^{n-1}}(X+1)^{p^{n-1}}$, we can equivalently reformulate (II) to
$$
F_n(X)\;\;\con_{p^n}\;\; F_{n-1}(X^p)\; .
\leqno (\mb{II}')
$$
\end{footnotesize}
\end{quote}

\subsection{Characteristic $p$}

In this section, we assume $\ch K = p$. 

\begin{Assumption}
\label{Ass3a}\rm
Suppose given $x,y\in T$ and $k\in [1,l-1]$ such that
\begin{itemize}
\item[(i)] $t^m\con_{ys} s$,
\item[(ii)] $x\; |\; a_j y^{j[p]}$ for all $j\in [1,l - 1]$, and
\item[(iii)] $x r\; |\; a_j y^{j[p]}$ for all $j\in [1,k - 1]$.
\end{itemize}
\end{Assumption}

Let $c := \gcd(x s^k,y^{l[p]} s^l)\in T$.

\begin{Lemma}
\label{Lem5a}
Given {\rm (\ref{Ass3a}),} we have $\; c\; |\; \mu_{s,K}(t^m)\;$.

\rm
We may decompose
\[
\begin{array}{rcl}
\mu_{s,K}(t^m)
& = & \mu_{s,K}(t^m) - \mu_{s,K}(s) \\
& = & (t^{ml} - s^l) + \Big(\sum_{j\in [1,k-1]} a_j (t^{mj} - s^j)\Big) + \Big(\sum_{j\in [k,l-1]} a_j (t^{mj} - s^j)\Big)\; . \\
\end{array}
\]
Now since $t^m\con_{ys} s$, we have $t^{mj}\con_{y^{j[p]} s^j} s^j$ for any $j\geq 1$.

In particular, $y^{l[p]} s^l\; |\; t^{ml} - s^l$.

Moreover, $x s^l\; |\;\sum_{j\in [1,k-1]} a_j (t^{mj} - s^j)$ by (\ref{Ass3a}.iii).

Finally, $x s^k\; |\; \sum_{j\in [k,l-1]} a_j (t^{mj} - s^j)$ by (\ref{Ass3a}.ii). 
\end{Lemma}

\begin{Proposition}
\label{Prop6a}
Given {\rm (\ref{Ass3a})}, we have $t^{-j} c \; |\; b_j$ if $j\not\con_m 0$ and $t^{-j} c \; |\; (b_j - a_{j/m})$ if $j\con_m 0$ for $j\in [1,lm - 1]$.

\rm
This follows using (\ref{Lem5a}), cf.\ (\ref{Prop6}).
\end{Proposition}

%% file: wsec2.tex
\section{A tower of purely ramified extensions}
\label{SecTow}

Suppose given a chain
\[
R_0\;\tm\; R_1\;\tm\; R_2\;\tm\;\cdots
\]
of finite purely ramified extensions $R_{i+1}|R_i$, with maximal ideal generated by $r_i\in R_i$, of residue characteristic $\ch R_i/r_i R_i = p$, with field of 
fractions $K_i = \fracfield R_i$, and of degree $[K_{i+1}:K_i] = p^\kappa = q$ for $i\geq 0$, where $\kappa\geq 1$ is an integer stipulated to be independent of 
$i$. We may and will suppose that $\Nrm_{K_{i+1}|K_i}(r_{i+1}) = r_i$ for $i\geq 0$. We write
\[
\mu_{r_i,K_0}(X)\; =\; X^{q^i} + \Big(\sum_{j\in [1,q^i - 1]} a_{i,j} X^j\Big) - r_0\;\in\; R_0[X]\; .
\]
For $j\geq 1$, we denote $v_q(j) := \max\{ \al\in\Z_{\geq 0}\; :\; j\con_{q^\al} 0 \}$. That is, $v_q(j)$ is the 
largest integer below $v_p(j)/\kappa$. We abbreviate $g := (q-2)/(q-1)$.

\begin{Assumption}
\label{Ass6_5}\rm
Suppose given $f\in R_0$ such that $r_i^{q-1} f\; |\; r_i^q - r_{i-1}$ for all $i\geq 0$. \linebreak If $\ch K_0 = 0$, then suppose $p\; |\; f\; |\; q$. 
If $\ch K_0 = p$, then suppose $r_0\; |\; f$.
\end{Assumption}

\begin{Proposition}
\label{Prop7}
Assume {\rm (\ref{Ass6_5}).} 
\begin{itemize}
\item[{\rm (i)}] We have $f^{i - v_q(j)}\; |\; a_{i,j}$ for $i\geq 1$ and $j\in [1,q^i-1]$. 
\item[{\rm (i$'$)}] If $j < q^i g$, then $f^{i - v_q(j)} r_0\; |\; a_{i,j}$.
\item[{\rm (ii)}] We have $a_{i,j} \con_{f^{i+1}} a_{i+\be,q^\be j}$ for $i\geq 1$, $j\in [1,q^i-1]$ and $\be\geq 1$. 
\item[{\rm (ii$'$)}] If $j < q^i g$, then $a_{i,j} \con_{f^{i+1}r_0} a_{i+\be,q^\be j}$ for $\be\geq 1$.
\end{itemize}

\rm
{\it Consider the case $\ch K_0 = 0$.} To prove (i, i$'$), we perform an induction on $i$, the assertion being true for $i = 1$ by 
(\ref{Ass6_5}). So suppose given $i\geq 2$ and the assertion to be true for $i - 1$. To apply (\ref{Prop6}), we let $R = R_0$, $r = r_0$, $S = R_{i-1}$, 
$s = r_{i-1}$, $T = R_i$ and $t = r_i$. Furthermore, we let $y = r_i^{q-1} f$, $x = f^{i-1}$ and $k = q^{i-1} - (q^{i-1} - 1)/(q - 1)$, so that (\ref{Ass3}) is 
satisfied by (\ref{Ass6_5}) and by the inductive assumption. We have $c = f^i r_i^{qk - 1}$.

Consider $j\in [1,q^i - 1]$. If $j\not\con_q 0$, then (\ref{Prop6}) gives 
\[
v_{r_i}(a_{i,j}/f^i) \;\geq\; qk - 1 - j\; , 
\]
whence $f^i$ divides $a_{i,j}$; $f^i$ strictly divides $a_{i,j}$, if $j < q^i g$ since $0 < (qk - 1) - q^i g = 1/(q - 1) < 1$.

If $j\con_q 0$, then (\ref{Prop6}) gives 
\[
v_{r_i}((a_{i,j} - a_{i-1,j/q})/f^i) \;\geq\; qk - 1 - j\; , 
\]
whence $f^i$ divides $a_{i,j} - a_{i-1,j/q}$; strictly, if $j < q^i g$. By induction, $f^{i - 1 - v_q(j/q)}$ divides $a_{i-1,j/q}$; strictly, if 
$j/q < q^{i-1}g$. But $a_{i-1,j/q} \con_{f^i} a_{i,j}$, and therefore $f^{i - v_q(j)}$ divides also $a_{i,j}$; strictly, if $j < q^i g$. This proves (i, i$'$).

The case $\be = 1$ of (ii, ii$'$) has been established in the course of the proof of (i, i$'$). The general case follows by induction.

{\it Consider the case $\ch K_0 = p$.} To prove (i, i$'$), we perform an induction on $i$, the assertion being true for $i = 1$ by (\ref{Ass6_5}). So suppose 
given $i\geq 2$ and the assertion to be true for $i - 1$. To apply (\ref{Prop6a}), we let $R = R_0$, $r = r_0$, $S = R_{i-1}$, $s = r_{i-1}$, $T = R_i$ and 
$t = r_i$. Furthermore, we let $y = r_i^{-1} f$, $x = r_i^{-1} f^i$ and $k = q^{i-1} - (q^{i-1} - 1)/(q - 1)$, so that (\ref{Ass3a}) is satisfied by 
(\ref{Ass6_5}) and by the inductive assumption. In fact, $x y^{-j[p]} = r_i^{j[p] - 1} f^{i - j[p]}$ divides $f^{i - 1 - v_q(j)}$ both if $j\not\con_p 0$ and if 
$j\con_p 0$; in the latter case we make use of the inequality $p^{\al - 1} (p-1) \geq \al + 1$ for $\al\geq 1$, which needs $p\geq 3$. We obtain 
$c = f^i r_i^{qk - 1}$. 

Using (\ref{Prop6a}) instead of (\ref{Prop6}), we may continue as in the former case to prove (i, i$'$), and, in the course of this proof, also (ii, ii$'$).
\end{Proposition}

%% file: wsec3.tex
\section{Galois descent of a divisibility}
\label{SecDes}

Let
\begin{center}
\begin{picture}(250,250)
\put(   0,   0){$S$}
\put(  60,   5){\vector(1,0){120}}
\put(  60,  15){\oval(20,20)[l]}
\put( 110,  25){$\scm G$}
\put( 200,   0){$\w S$}
\put(  15,  60){\vector(0,1){120}}
\put(   5,  60){\oval(20,20)[b]}
\put(  30, 110){$\scm m$}
\put( 215,  60){\vector(0,1){120}}
\put( 205,  60){\oval(20,20)[b]}
\put( 230, 110){$\scm m$}
\put(   0, 200){$T$}
\put(  60, 205){\vector(1,0){120}}
\put(  60, 215){\oval(20,20)[l]}
\put( 110, 225){$\scm G$}
\put( 200, 200){$\w T$}
\end{picture}
\end{center}
be a commutative diagram of finite, purely ramified extensions of discrete valuation rings. Let $s\in S$, $t\in T$, $\w s\in \w S$ and $\w t\in \w T$ generate 
the respective maximal ideals. Let $L = \fracfield S$, $M = \fracfield T$, $\w L = \fracfield \w S$ and $\w M = \fracfield \w T$ denote the respective field of 
fractions. We assume the extensions $M|L$ and $\w L|L$ to be linearly disjoint and $\w M$ to be the composite of $M$ and $\w L$. Thus $m := [M:L] = [\w M:\w L]$ 
and $[\w L : L] = [\w M : M]$. We assume $\w L|L$ to be galois and identify $G := \Gal(\w L|L) = \Gal(\w M|M)$ via restriction. We may and will assume that 
$s = \Nrm_{\w L|L}(\w s)$, and that $t = \Nrm_{\w M|M}(\w t)$.

\begin{Lemma}
\label{Lem9}
In $\;\w T\,$, the element $\; 1 - \w t^m/\w s\;$ divides $\; 1 - t^m/s\;$.

\rm
Let $\w d = 1 - \w t^m/\w s$, so that $\w t^m = \w s (1 - \w d)$. We conclude
\[
\begin{array}{rll}
t^m
& = & \Nrm_{\w M|M}(\w t^m) \\
& = & \Nrm_{\w L|L}(\w s)\cdot\prod_{\sigma\in G} (1 - \w d^\sigma) \\
& \con_{s\w d} & s\; . \\
\end{array}
\]
\end{Lemma}

%% file: wsec4.tex
\section{Cyclotomic number fields}

\subsection{Coefficient valuation bounds}
\label{SubSecMipo}

For $n \geq 1$, we let $\zeta_{p^n}$ be a primitive $p^n$th root of unity over $\Q$. We make choices in such a manner that $\zeta_{p^n}^p = \zeta_{p^{n-1}}$
for $n\geq 2$. We denote $\th_n = \zeta_{p^n} - 1$ and $F_n = \Q(\zeta_{p^n})$. Let 
$E_n = \mb{Fix}_{C_{p-1}} F_n$, so $[E_n:\Q] = p^{n-1}$. Let 
\[
\pi_n \; =\; \Nrm_{F_n|E_n}(\th_n) \; =\; \prod_{j\in [1,\, p-1]} (\zeta_{p^n}^{j^{p^{n-1}}} - 1)\; .
\]
The minimal polynomial $\mu_{\th_n,F_{n-1}}(X) = (X+1)^p - \th_{n-1} - 1$ shows that $\Nrm_{F_n|F_{n-1}}(\th_n) = \th_{n-1}$, hence also 
$\Nrm_{E_n|E_{n-1}}(\pi_n) = \pi_{n-1}$. Note that $\pi_1 = p$ and $E_1 = \Q$.

Let $\Ol$ be the integral closure of $\Z_{(p)}$ in $E_n$. Since $\Nrm_{E_n|\sQ}(\pi_n) = \pi_1 = p$, we have $\Z_{(p)}/p\Z_{(p)}\lraiso \Ol/\pi_n\Ol$. In 
particular, the ideal $\pi_n\Ol$ in $\Ol$ is prime. Now $\pi_n^{p^{n-1}}\Ol = p\Ol$, since 
$\pi_n^{p^{n-1}}/p = \pi_n^{p^{n-1}}/\Nrm_{E_n|\sQ}(\pi_n)\in \Z_{(p)}[\th_n]^\ast\cap E_n = \Ol^\ast$. Thus $\Ol$ is a discrete valuation ring, purely 
ramified of degree $p^{n-1}$ over $\Z_{(p)}$, and so $\Ol = \Z_{(p)}[\pi_n]$ {\bf \cite[\rm I.\S 7, prop.\ 18]{Se62}}. In particular, $E_n = \Q(\pi_n)$. 

\begin{Lemma}
\label{Lem10}
For all $n\geq 2$, we have $\;\pi_n^p\;\con_{\pi_n^{p-1} p}\;\pi_{n-1}\;$.

\rm
First of all, $\th_n^p\con_{\th_n p}\th_{n-1}$ since $(X - 1)^p - (X^p - 1)$ is divisible by $p(X - 1)$ in $\Z[X]$. Letting $\w T = \Z_{(p)}[\th_n]$ and 
$(\w t,\w s,t,s) = (\th_n,\th_{n-1},\pi_n,\pi_{n-1})$, (\ref{Lem9}) shows that $1 - \th_n^p/\th_{n-1}$ divides $1 - \pi_n^p/\pi_{n-1}$. Therefore, 
$\th_n p\, \th_{n-1}^{-1}\pi_{n-1} \; |\; \pi_{n-1} - \pi_n^p$.
\end{Lemma}

Now suppose given $m\geq 1$. To apply (\ref{Prop7}), we let $f = q = p$, $R_i = \Z_{(p)}[\pi_{m+i}]$ and $r_i = \pi_{m+i}$ for $i\geq 0$. We keep the
notation
\[
\mu_{\pi_{m+i},\,E_m}(X)\; =\;\mu_{r_i,\,K_0}(X)\; =\; X^{p^i} + \Big(\sum_{j\in [1,p^i - 1]} a_{i,j} X^j\Big) - \pi_m \;\in\; R_0[X]
\; = \; \Z_{(p)}[\pi_m][X]\; .
\]

\begin{Theorem}
\label{Th11}
\Absit
\begin{itemize}
\item[{\rm (i)}] We have $p^i\; |\; j a_{i,j}$ for $i\geq 1$ and $j\in [1,p^i-1]$. 
\item[{\rm (i$'$)}] If $j < p^i (p-2)/(p-1)$, then $p^i\pi_m\; |\; j a_{i,j}$.
\item[{\rm (ii)}] We have $a_{i,j} \con_{p^{i+1}} a_{i+\be,p^\be j}$ for $i\geq 1$, $j\in [1,p^i-1]$ and $\be\geq 1$. 
\item[{\rm (ii$'$)}] If $j < p^i (p-2)/(p-1)$, then $a_{i,j} \con_{p^{i+1}\pi_m} a_{i+\be,p^\be j}$.
\end{itemize}

\rm
Assumption (\ref{Ass6_5}) is fulfilled by virtue of (\ref{Lem10}), whence the assertions follow by (\ref{Prop7}).
\end{Theorem}

\begin{quote}
 \begin{footnotesize}
 \begin{Example}
 \label{ExEis12}\rm
 For $p = 5$, $m = 1$ and $i = 2$, we have 
 \[
 \begin{array}{rl}
 \mu_{\pi_3,\sQ}(X)\;\; = & X^{25} - 4\cdot 5^2 X^{24} + 182\cdot 5^2 X^{23} - 8\cdot 5^6 X^{22} + 92823\cdot 5^2 X^{21} \\
                      & -\; 6175454\cdot 5 X^{20} + 12194014 \cdot 5^2 X^{19} - 18252879\cdot 5^3 X^{18} \\
                      & +\; 4197451\cdot 5^5 X^{17} - 466901494\cdot 5^3 X^{16} + 8064511079\cdot 5^2 X^{15}\\ 
                      & -\; 4323587013\cdot 5^3 X^{14} + 1791452496\cdot 5^4 X^{13} - 113846228\cdot 5^6 X^{12} \\
                      & +\; 685227294\cdot 5^5 X^{11} - 15357724251\cdot 5^3 X^{10} + 2002848591\cdot 5^4 X^9 \\ 
                      & -\; 4603857997\cdot 5^3 X^8 + 287207871\cdot 5^4 X^7 - 291561379\cdot 5^3 X^6 \\
                      & +\; 185467152\cdot 5^2 X^5 - 2832523 \cdot 5^3 X^4 + 121494 \cdot 5^3 X^3 - 514\cdot 5^4 X^2\\ 
                      & +\; 4\cdot 5^4 X - 5\; . \\
 \end{array}  
 \]
 Now $v_5(a_{3,22}) = 6\neq 5 = v_5(a_{4,5\cdot 22})$, so the valuations of the coefficients considered in (\ref{Th11}.ii) differ in general. This, however, 
 does not contradict the assertion $a_{3,22}\con_{5^4} a_{4,\, 5\cdot 22}$ from loc.\ cit.
 \end{Example}

 \end{footnotesize}
\end{quote}

\subsection{A different proof of (\ref{Th11}.\,i,\,i$'$) and some exact valuations}
\label{SubSecAlt}

Let $m\geq 1$ and $i\geq 0$. We denote $R_i = \Z_{(p)}[\pi_{m+i}]$, $r_i = \pi_{m+i}$, $K_i = \fracfield R_i$, \linebreak $\w R_i = \Z_{(p)}[\th_{m+i}]$ and 
$\w r_i = \th_{m+i}$. Denoting by $\Dfk$ the respective different {\bf\cite[\rm III.\S3]{Se62}}, we have $\Dfk_{\w R_i|\w R_0} = (p^i)$ and
$\Dfk_{\w R_i|R_i} = (\w r_i^{p-2})$ {\bf\cite[\rm III.\S3, prop.\ 13]{Se62}}, whence
$$
\Dfk_{R_i|R_0} \= \left(\mu'_{r_i,K_0}(r_i)\right) \= \Dfk_{\w R_i|\w R_0}\Dfk_{\w R_0|R_0}\Dfk_{\w R_i|R_i}^{-1} 
\= \left(p^i r_i^{p^i - 1 - (p^i - 1)/(p - 1)}\right)\; ,
\leqno (\ast)
$$
cf.\ {\bf\cite[\rm III.\S3, cor.\ 2]{Se62}}. Therefore, $p^i r_i^{p^i - 1 - (p^i - 1)/(p - 1)}$ divides $j a_{i,j} r_i^{j-1}$ for $j\in [1,p^i - 1]$, and
(\ref{Th11}.\,i,\,i$'$) follow.

Moreover, since only for $j = p^i - (p^i - 1)/(p - 1)$ the valuations at $r_i$ of $p^i r_i^{p^i - 1 - (p^i - 1)/(p - 1)}$ and $j a_{i,j} r_i^{j-1}$
are congruent modulo $p^i$, we conclude by ($\ast$) that they are equal, i.e.\ that
\[
a_{i,p^i - (p^i - 1)/(p - 1)}\,/\, p^i\;\in\; R_0^\ast\; .
\]

\begin{Corollary}
\label{CorDiff}
We have
\[
a_{i,p^i - (p^i - p^\be)/(p - 1)}\,/\, p^{i-\be}\;\in\; R_0^\ast\= \Z_{(p)}[\pi_m]^\ast\Icm\mb{for $\;\;\be\in [0,i-1]$}\; .
\]

\rm
This follows by (\ref{Th11}.ii) from what we have just said.
\end{Corollary}

E.g.\ in (\ref{ExEis12}), $5^1$ exactly divides $a_{2,25-5} = a_{2,20}$, and $5^2$ exactly divides $a_{2,25-5-1} = a_{2,19}$.

\subsection{Some traces}
\label{SubsecTr}
Let $\mufat_{p-1}$ denote the group of $(p-1)$st roots of unity in $\Q_p.$
We choose a primitive $(p-1)$st root of unity $\zeta_{p-1}\in\mufat_{p-1}$ and may thus view $\Q(\zeta_{p-1})\tm\Q_p$ as a 
subfield. Note that $[\Q(\zeta_{p-1}) : \Q] = \phi(p-1)$, where $\phi$ denotes Euler's function. The restriction of the valuation $v_p$ at $p$ on $\Q_p$ to 
$\Q(\zeta_p)$, is a prolongation of the valuation $v_p$ on $\Q$ to $\Q(\zeta_{p-1})$ (there are $\phi(p-1)$ such prolongations).
 
\begin{Proposition}
\label{PropFT2}
For $n\geq 1$, we have 
\[
\Tr_{E_n|\sQ}(\pi_n)\; =\; p^n s_n - p^{n-1} s_{n-1}\; , 
\]
where
\[
s_n \; :=\;  \frac{1}{p-1}\sum_{H\,\tm\,\mufat_{p-1}}(-1)^{\# H}{\txm\Big\{v_p\big(\sum_{\xi\in H}\xi\big)\ge n\Big\}} \quad\text{for}\quad n\ge 0\; .
\]
We have $s_0 = 0$, and $s_n\in\Z$ for $n\geq 0$. The sequence $(s_n)_n$ becomes stationary at some minimally chosen $N_0(p)$. We have 
\[
N_0(p)\;\leq\;
N(p)\; :=\; \max_{H\,\tm\,\mufat_{p-1}}\Big\{v_p\big({\txm\sum_{\xi\in H}\xi}\big)\!:\,{\txm\sum_{\xi\in H} \xi} \neq 0\Big\}+1\; .
\] 

\rm
An upper estimate for $N(p)$, hence for $N_0(p)$, is given in (\ref{LemUB}).

{\it Proof of {\rm (\ref{PropFT2}).}}
For $j\in [1,p-1]$ the $p$\,-adic limits
$$
\xi(j)\; :=\; \lim_{n\to\infty}\, j^{p^n}
$$
exist since  $j^{p^{n-1}}\con_{p^n} j^{p^n}$. They are distinct since $\xi(j)\con_{p}j$, and, thus, form the group $\mufat_{p-1} = \{\xi(j)\; |\; j\in [1,p-1]\}$.
Using the formula
$$
\Tr_{F_n|\sQ}(\zeta_{p^n}^m)\= p^n\big\{v_p(m)\ge n\big\} - p^{n-1}\big\{v_p(m)\ge n-1\big\}
$$
and the fact that $j^{p^{n-1}}\con_{p^n}\xi(j)$, we obtain 
\[
\begin{array}{rcl}
\Tr_{F_n|\sQ}(\pi_n)
& = & \Tr_{F_n|\sQ}\bigg(\prod_{j\in [1,\, p-1]}\Big(1 - \zeta_{p^n}^{j^{p^{n-1}}}\Big)\bigg) \\
& = & \sumd{J\,\tm\, [1,\, p-1]} (-1)^{\# J}\,\Tr_{F_n|\sQ}\Big(\zeta_{p^n}^{\sum_{j\in J} j^{p^{n-1}}}\Big) \\
& = & \sumd{J\,\tm\, [1,\, p-1]} (-1)^{\# J}\bigg(p^n \left\{ v_p{\txm\big(\sum_{j\in J}\xi(j)\big)}\geq n \right\} \\
&   & \hspace*{23.7mm} -\; p^{n-1}{\txm\Big\{v_p\big(\sum_{j\in J} \xi(j)\big)\geq n-1\Big\}} \bigg) \\
& = & (p-1)(p^n s_n - p^{n-1} s_{n-1})\; , \\
\end{array}
\]
whence 
\[
\Tr_{E_n|\sQ}(\pi_n) \= p^n s_n - p^{n-1} s_{n-1}\; .
\]
Now $s_0 = 0\in\Z$ by the binomial formula. Therefore, by induction, we conclude from $p^n s_n - p^{n-1} s_{n-1}\in\Z$ that $p^n s_n\in\Z$. 
Since $(p - 1)s_n\in\Z$, too, we obtain $s_n\in\Z$.

As soon as $n\ge N(p)$, the conditions $v_p(\sum_{\xi\in H} \xi) \geq n$ and  $v_p(\sum_{\xi\in H} \xi) = +\infty$ on $H\tm\mufat_{p-1}$ become equivalent, 
and we obtain
$$
s_n\= \frac{1}{p-1}\sum_{H\tm\mufat_{p-1}}(-1)^{\#H}\big\{{\txm\sum_{\xi\in H}\xi = 0}\big\}\,,
$$
which is independent of $n$. Thus $N_0(p)\le N(p)$. 
\end{Proposition}

\begin{Lemma}
\label{LemFT3}
We have $s_1 = 1$. In particular,
$\;
\Tr_{E_2|\sQ}(\pi_2)\;\con_{p^2}\; -p \; . 
$

\rm
Since $\Tr_{E_1|\sQ}(\pi_1) = \Tr_{\sQ|\sQ}(p) = p$, and since $s_0 = 0$, we have $s_1 = 1$ by (\ref{PropFT2}). The congruence for 
$\Tr_{E_2|\sQ}(\pi_2)$ follows again by (\ref{PropFT2}).
\end{Lemma}

\begin{Corollary}
\label{CorEis13}
We have
\[
\mu_{\pi_n,\sQ}(X)\; \con_{p^2} X^{p^{n-1}} + pX^{(p-1)p^{n-2}} - p\;
\]
for $n\geq 2$.

\rm
By dint of (\ref{LemFT3}), this ensues from (\ref{Th11}.\,i$'$,\,ii).
\end{Corollary}

\begin{quote}
\begin{footnotesize}
\begin{Example}
\label{ExEis7_5}\rm
The last $n$ for which we list $s_n$ equals $N(p)$, except if there is a question mark in the next column. The table was calculated using Pascal ($p\leq 53$)
and Magma ($p\geq 59$). In the last column, we list the upper bound for $N(p)$ calculated below (\ref{LemUB}).

\begin{scriptsize}
\hspace*{-10mm}
\begin{tabular}{r||r|r|r|r|r|r|r|r||c}
 $s_n$   & $n = 0$             & $1$               & $2$               &  $3$ & $4$ & $5$  & $6$ & $7$ & \shortstack{upper bound\\for $N(p)$}\\\hline\hline
 $p = 3$ & $0$                 & $1$               &                   &                &                &                &           &           &       \\\hline
 $5$     & $0$                 & $1$               &                   &                &                &                &           &           & $1$   \\\hline
 $7$     & $0$                 & $1$               &                   &                &                &                &           &           & $1$   \\\hline
 $11$    & $0$                 & $1$               & $3$               &                &                &                &           &           & $3$   \\\hline
 $13$    & $0$                 & $1$               & $3$               &                &                &                &           &           & $3$   \\\hline
 $17$    & $0$                 & $1$               & $8$               & $16$           &                &                &           &           & $5$   \\\hline
 $19$    & $0$                 & $1$               & $10$              & $12$           &                &                &           &           & $4$   \\\hline
 $23$    & $0$                 & $1$               & $33$              & $89$           & $93$           &                &           &           & $7$   \\\hline
 $29$    & $0$                 & $1$               & $377$             & $571$          & $567$          &                &           &           & $8$   \\\hline
 $31$    & $0$                 & $1$               & $315$             & $271$          & $259$          &                &           &           & $6$   \\\hline
 $37$    & $0$                 & $1$               & $107$             & $940$          & $1296$         &                &           &           & $9$   \\\hline
 $41$    & $0$                 & $1$               & $6621$            & $51693$        & $18286$        & $20186$        & $20250$   &           & $12$  \\\hline
 $43$    & $0$                 & $1$               & $1707$            & $4767$         & $6921$         & $6665$         &           &           & $9$   \\\hline
 $47$    & $0$                 & $1$               & $2250$            & $272242$       & $173355$       & $181481$       & $182361$  &           & $16$  \\\hline
 $53$    & $0$                 & $1$               & $71201$           & $363798$       & 1520045        & $1350049$      & $1292229$ & $1289925$ & $18$  \\\hline
 $59$    & $0$                 & $1$               & $1276$            & ?              &                &                &           &           & $21$  \\\hline
 $61$    & $0$                 & $1$               & $2516$            & ?              &                &                &           &           & $12$  \\\hline
 $67$    & $0$                 & $1$               & $407186$          & ?              &                &                &           &           & $15$  \\\hline
 $71$    & $0$                 & $1$               & $5816605$         & ?              &                &                &           &           & $18$  \\\hline
 $73$    & $0$                 & $1$               & $8370710$         & ?              &                &                &           &           & $18$  \\\hline
 $79$    & $0$                 & $1$               & $169135$          & ?              &                &                &           &           & $18$  \\\hline
 $83$    & $0$                 & $1$               & $632598$          & ?              &                &                &           &           & $30$  \\\hline
 $89$    & $0$                 & $1$               & $26445104$        & ?              &                &                &           &           & $30$  \\\hline
 $97$    & $0$                 & $1$               & $282789$          & ?              &                &                &           &           & $24$  \\\hline
 $101$   & $0$                 & $1$               & $25062002$        & ?              &                &                &           &           & $31$  \\\hline
 $103$   & $0$                 & $1$               & $56744199$        & ?              &                &                &           &           & $25$  \\\hline
 $107$   & $0$                 & $1$               & $1181268305$      & ?              &                &                &           &           & $40$  \\\hline
 $109$   & $0$                 & $1$               & $91281629$        & ?              &                &                &           &           & $28$  \\\hline
 $113$   & $0$                 & $1$               & $117774911422$    & ?              &                &                &           &           & $37$  \\\hline
 $127$   & $0$                 & $1$               & $6905447$         & ?              &                &                &           &           & $28$  \\\hline
 $131$   & $0$                 & $1$               & $2988330952791$   & ?              &                &                &           &           & $37$  \\\hline
 $137$   & $0$                 & $1$               & $1409600547$      & ?              &                &                &           &           & $50$  \\\hline
 $139$   & $0$                 & $1$               & $3519937121$      & ?              &                &                &           &           & $34$  \\\hline
 $149$   & $0$                 & $1$               & $25026940499$     & ?              &                &                &           &           & $56$  \\\hline
 $151$   & $0$                 & $1$               & $164670499159$    & ?              &                &                &           &           & $31$  \\\hline
 $157$   & $0$                 & $1$               & $51594129045351$  & ?              &                &                &           &           & $38$  \\\hline
 $163$   & $0$                 & $1$               & $288966887341$    & ?              &                &                &           &           & $42$  \\\hline
 $167$   & $0$                 & $1$               & $1205890070471$   & ?              &                &                &           &           & $64$  \\\hline
 $173$   & $0$                 & $1$               & $17802886165762$  & ?              &                &                &           &           & $66$  \\\hline
 $179$   & $0$                 & $1$               & $1311887715966$   & ?              &                &                &           &           & $69$  \\\hline
 $181$   & $0$                 & $1$               & $128390222739$    & ?              &                &                &           &           & $38$  \\\hline
 $191$   & $0$                 & $1$               & $233425263577158$ & ?              &                &                &           &           & $57$  \\\hline
 $193$   & $0$                 & $1$               & $306518196952028$ & ?              &                &                &           &           & $51$  \\\hline
 $197$   & $0$                 & $1$               & $347929949728221$ & ?              &                &                &           &           & $66$  \\\hline
 $199$   & $0$                 & $1$               & $9314622093145$   & ?              &                &                &           &           & $48$  \\\hline
 $211$   & $0$                 & $1$               & $12532938009082$  & ?              &                &                &           &           & $39$  \\\hline
\end{tabular}
\end{scriptsize}%

So for example if $p = 31$, then $\Tr_{\sQ(\pi_3)|\sQ}(\pi_3) = 271\cdot 31^3 - 315\cdot 31^2$, whereas 
$\Tr_{\sQ(\pi_7)|\sQ}(\pi_7) = 259\cdot 31^7 - 259\cdot 31^6$. Moreover, $N_0(31) = N(31) = 4 \leq 6$.
\end{Example}

\begin{Remark}
\label{RemFT4}\rm
Vanishing (resp.\ vanishing modulo a prime) of sums of roots of unity has been studied extensively. See e.g.\ 
{\bf\cite{DZ02}}, {\bf\cite{LL00}}, where also further references may be found.
\end{Remark}

\begin{Remark}
\label{RemFT5}\rm
Neither do we know whether $s_n\geq 0$ nor whether $\Tr_{E_n|\sQ}(\pi_n)\geq 0$ always hold. Moreover, we do not know a prime $p$ for which 
$N_0(p) < N(p)$.
\end{Remark}

\begin{Remark}
\label{RemFT6}\rm
We calculated some further traces appearing in (\ref{Th11}), using Maple and Magma. 

For $p = 3$, $n\in [2,10]$, we have $\Tr_{E_n|E_{n-1}}(\pi_n) = 3\cdot 2$. 

For $p = 5$, $n\in [2,6]$, we have $\Tr_{E_n|E_{n-1}}(\pi_n) = 5\cdot 4$.

For $p = 7$, $n\in [2,5]$, we have $\Tr_{E_n|E_{n-1}}(\pi_n) = 7\cdot 6$.
  
For $p = 11$, we have $\Tr_{E_2|E_1}(\pi_2) = 11\cdot 32$, whereas 
\[
\begin{array}{cl}
    & \hspace*{-7mm}\Tr_{E_3|E_2}(\pi_3) \\
  = & 22\cdot (15 + \zeta^2 + 2\zeta^3 - \zeta^5 + \zeta^6 - 2\zeta^8 - \zeta^9 + 2\zeta^{14} - \zeta^{16} + \zeta^{18} - \zeta^{20} - 2\zeta^{24} \\
    & + 2\zeta^{25} - 2\zeta^{26} - \zeta^{27} - \zeta^{31} + 2\zeta^{36} - \zeta^{38} + \zeta^{41} - \zeta^{42} - 2\zeta^{43} + 2\zeta^{47} - 3\zeta^{49} \\
    & - \zeta^{53} + \zeta^{54} + 2\zeta^{58} - \zeta^{60} - \zeta^{64} + \zeta^{67}  + 2\zeta^{69} - \zeta^{71} - 2\zeta^{72} - \zeta^{75} - 2\zeta^{78} \\
    & + 3\zeta^{80} - \zeta^{82} - \zeta^{86} + 2\zeta^{91} - \zeta^{93} - 2\zeta^{95} - 3\zeta^{97} + 2\zeta^{102} + \zeta^{103} - \zeta^{104} - \zeta^{108}) \\  
  = & 22\cdot 2014455354550939310427^{-1}\cdot (34333871352527722810654 \\
    & + 1360272405267541318242502 \pi - 31857841148164445311437042 \pi^2 \\
    & + 135733708409855976059658636 \pi^3 - 83763613130017142371566453 \pi^4 \\
    & + 20444806599344408104299252\pi^5 - 2296364631211442632168932 \pi^6 \\
    & + 117743741083866218812293 \pi^7 - 2797258465425206085093 \pi^8 \\
    & + 27868038642441136108 \pi^9 - 79170513243924842\pi^{10})\; , \\
\end{array}
\]
where $\zeta := \zeta_{11^2}$ and $\pi := \pi_2$. 
\end{Remark}
\end{footnotesize}
\end{quote}

\subsection{An upper bound for $N(p)$}
\label{SubsecNp}

We view $\Q(\zeta_{p-1})$ as a subfield of $\Q_p$, and now, in addition, as a subfield of $\C$. Since complex conjugation commutes with the operation of 
$\Gal(\Q(\zeta_{p-1})|\Q)$, we have $|\Nrm_{\sQ(\zeta_{p-1})|\sQ}(x)| = |x|^{\phi(p-1)}$ for $x\in\Q(\zeta_{p-1})$.

We abbreviate $\Sigma(H) := \sum_{\xi\in H} \xi$ for $H\tm\mufat_{p-1}$. Since $|\Sigma(H)| \leq p-1$, we have \linebreak[4]
$|\Nrm_{\sQ(\zeta_{p-1})|\sQ}(\Sigma(H))|\leq (p-1)^{\phi(p-1)}$. Hence, if $\Sigma(H)\neq 0$, then 
\[
v_p(\Sigma(H)) \;\leq\; v_p(\Nrm_{\sQ(\zeta_{p-1})|\sQ}(\Sigma(H))) \; <\; \phi(p-1)\; , 
\]
and therefore $N(p) \leq \phi(p-1)$. We shall ameliorate this bound by a logarithmic term.

\begin{Proposition}
\label{LemUB}
We have
\[
N(p)\;\leq\; \phi(p-1)\left(1 - \frac{\log\pi}{\log p}\right) + 1
\]
for $p\geq 5$. 

\rm
It suffices to show that $|\Sigma(H)|\leq p/\pi$ for $H\tm\mufat_{p-1}$. We will actually show that
\[
\max_{H\tm\mufat_{p-1}} |\Sigma(H)|\; =\; \frac{1}{\sin\frac{\pi}{p-1}}\; ,
\]
from which this inequality follows using $\sin x\geq x - x^3/6$ and $p\geq 5$.

Choose $H\tm\mufat_{p-1}$ such that $|\Sigma(H)|$ is maximal. Since $p - 1$ is even, the $(p-1)$st roots of unity fall into pairs $(\eta,-\eta)$.
The summands of $\Sigma(H)$ contain exactly one element of each such pair, since $|\Sigma(H) + \eta|^2 + |\Sigma(H) - \eta|^2 = 2|\Sigma(H)|^2 + 2$
shows that at least one of the inequalities $|\Sigma(H) + \eta|\leq |\Sigma(H)|$ and $|\Sigma(H) - \eta|\leq |\Sigma(H)|$ fails.

By maximality, replacing a summand $\eta$ by $-\eta$ in $\Sigma(H)$ does not increase the value of $|\Sigma(H)|$, whence
\[
|\Sigma(H)|^2\;\geq\; |\Sigma(H) - 2\eta|^2 = |\Sigma(H)|^2 - 4\,\mb{Re}(\eta\cdot\ol{\Sigma(H)}) + 4\; ,
\]
and thus
\[
\mb{Re}(\eta\cdot\ol{\Sigma(H)}) \;\geq 1\; >\; 0\; .
\]
Therefore, the $(p-1)/2$ summands of $\Sigma(H)$ lie in one half-plane, whence the value of $|\Sigma(H)|$.
\end{Proposition}

%% file: wsec5.tex
\section{Cyclotomic function fields, after Carlitz and Hayes}

\subsection{Notation and basic facts}

\begin{quote}
\begin{footnotesize}
We shall give a brief review while fixing notation. 
\end{footnotesize}
\end{quote}

Let $\rho\geq 1$ and $r := p^\rho$. Write $\Zl := \Fu{r}[Y]$ and $\Ql := \Fu{r}(Y)$, where $Y$ is an independent variable. 
We fix an algebraic closure $\b\Ql$ of $\Ql$. The {\it Carlitz module structure} on $\b\Ql$ is defined by the $\Fu{r}$-algebra homomorphism given on the 
generator $Y$ as
\[
\begin{array}{rcl}
\Zl & \lra     & \End_{\Ql}\b\Ql \\
Y   & \lramaps & \Big(\xi\;\lramaps\; \xi^Y \; :=\; Y\xi + \xi^r \,\Big)\; . \\
\end{array}
\]
We write the module product of $\xi\in\b\Ql$ with $e\in\Zl$ as $\xi^e$. For each $e\in\Zl$, there exists
a unique polynomial $P_e(X)\in\Zl[X]$ that satisfies $P_e(\xi) = \xi^e$ for all $\xi\in\b\Ql$. In fact, $P_1(X) = X$, $P_Y(X) = YX + X^r$, and 
$P_{Y^{i+1}} = Y P_{Y^i}(X) +  P_{Y^i}(X^r)$ for $i\geq 1$. For a general $e\in\Zl$, the polynomial $P_e(Y)$ is given by the according linear combination of these.

Note that $P_e(0) = 0$, and that $P'_e(X) = e$, whence $P_e(X)$ is separable, i.e.\ it decomposes as a product of distinct linear factors in $\b\Ql[X]$. Let 
\[
\lafat_e \; =\; \mb{ann}_e \b\Ql\; =\; \{\xi\in\b\Ql\; :\; \xi^e = 0 \} \;\tm\;\b\Ql 
\]
be the annihilator submodule. Separability of $P_e(X)$ shows that $\#\lafat_e = \deg P_e(X) = r^{\deg e}$. Given a $\Ql$-linear automorphism $\sigma$ of 
$\b\Ql$, we have $(\xi^e)^\sigma = P_e(\xi)^\sigma = P_e(\xi^\sigma) = (\xi^\sigma)^e$. In particular, $\lafat_e$ is stable under $\sigma$. Therefore, 
$\Ql(\lafat_e)$ is a Galois extension of $\Ql$.

Since $\#\mb{ann}_{\w e}\lafat_e = \#\lafat_{\w e} = r^{\deg\w e}$ for $\w e\; |\; e$, we have $\lafat_e\iso \Zl/e$ as $\Zl$-modules. It is not possible, 
however, to distinguish a particular isomorphism.

We shall restrict ourselves to prime powers now.
We fix a monic irreducible polynomial $f = f(Y) \in\Zl$ and write $q := r^{\deg f}$.
For $n\geq 1$, we let $\theta_n$ be a $\Zl$-linear generator of $\lafat_{f^n}$. We make our choices in such a manner that 
$\theta_{n+1}^f = \theta_n$ for $n\geq 1$. Note that $\Zl[\lafat_{f^n}] = \Zl[\theta_n]$ since the elements of $\lafat_{f^n}$ are polynomial expressions in 
$\theta_n$. 

Suppose given two roots $\xi, \w\xi\in\b\Ql$ of 
\[
\Psi_{\! f^n}(X) \; := \; P_{\! f^n}(X)/P_{\! f^{n-1}}(X)\;\in\;\Zl[X]\; , 
\]
i.e.\ $\xi, \w\xi\in\lafat_{f^n}\ohne\lafat_{f^{n-1}}$. Since 
$\xi$ is a $\Zl$-linear generator of $\lafat_{f^n}$, there is an $e\in\Zl$ such that $\w\xi = \xi^e$. Since $\xi^e/\xi = P_e(X)/X|_{X = \xi}\in\Zl[\theta_n]$, 
$\w\xi$ is a multiple of $\xi$ in $\Zl[\theta_n]$. Reversing the argument, we see that $\w\xi$ is in fact a unit multiple of $\xi$ in $\Zl[\theta_n]$.

\begin{Lemma}
\label{LemCar2}
The polynomial $\Psi_{\! f^n}(X)$ is irreducible. 

\rm
We have $\Psi_{\! f^n}(0) = \left.\fracd{P_{\! f^n}(X)/X}{P_{\! f^{n-1}}(X)/X}\right|_{X = 0} = f$\vspace{1mm}. We decompose 
$\Psi_{\! f^n}(X) = \prod_{i\in [1,k]} F_i(X)$ in its distinct monic irreducible factors $F_i(X)\in\Zl[X]$. One of the constant terms, say $F_j(0)$, is thus a 
unit multiple of $f$ in $\Zl$, while the other constant terms are units. Thus, all roots of $F_j(X)$ in $\Ql[\theta_n]$ are non-units in $\Zl[\theta_n]$, 
and the remaining roots of $\Psi_{\! f^n}(X)$ are units. But all roots of $\Psi_{\! f^n}(X)$ are unit multiples of each other. We conclude that 
$\Psi_{\! f^n}(X) = F_j(X)$ is irreducible.
\end{Lemma}

By (\ref{LemCar2}), $\Psi_{\! f^n}(X)$ is the minimal polynomial of $\theta_n$ over $\Ql$. In particular, $[\Ql(\theta_n) : \Ql] = q^{n-1}(q - 1)$, and so 
\[
\Zl[\theta_n]\theta_n^{(q-1)q^{n-1}}\; =\;\Zl[\theta_n]\Nrm_{\Ql(\theta_n)|\Ql}(\theta_n) \; =\; \Zl[\theta_n]f\; . 
\]
In particular, $\Zl_{(f)}[\theta_n]$ is a discrete valuation ring with maximal ideal generated by $\theta_n$, purely ramified of index $q^{n-1}(q-1)$ over 
$\Zl_{(f)}$, cf.\ {\bf \cite[\rm I.\S 7, prop.\ 18]{Se62}}. There is a group isomorphism
\[
\begin{array}{rcl}
(\Zl/f^n)^\ast & \lraiso  & \Gal(\Ql(\theta_n)|\Ql) \\
e              & \lramaps & (\theta_n\lramaps \theta_n^e)\; , \\
\end{array}
\]
well defined since $\theta_n^e$ is a root of $\Psi_{\! f^n}(X)$, too; injective since $\theta_n$ generates $\lafat_{f^n}$ over $\Zl$; and surjective by 
cardinality.

Note that the Galois operation on $\Ql(\theta_n)$ corresponding to $e\in (\Zl/f^n)^\ast$ coincides with the module operation of $e$ on the element $\theta_n$, 
but not everywhere. For instance, if $f\neq Y$, then the Galois operation corresponding to $Y$ sends $1$ to $1$, whereas the module operation of $Y$ sends 
$1$ to $Y + 1$. 

The discriminant of $\Zl[\theta_n]$ over $\Zl$ is given by 
$
\Delta_{\Zl[\theta_n]|\Zl} 
 =  \Nrm_{\Ql(\theta_n)|\Ql}(\Psi'_{\! f^n}(\theta_n)) 
$\linebreak[4]
$
 =  \Nrm_{\Ql(\theta_n)|\Ql}\left(P'_{\! f^n}(\theta_n)/P_{\! f^{n-1}}(\theta_n)\right) 
 =  \Nrm_{\Ql(\theta_n)|\Ql}\left(f^n/\theta_1\right) 
 =  f^{q^{n-1}(nq - n - 1)}\; . 
$

\begin{Lemma}
\label{LemCar5}
The ring $\Zl[\theta_n]$ is the integral closure of $\Zl$ in $\Ql(\theta_n)$. 

\rm
Let $e\in\Zl$ be a monic irreducible polynomial different from $f$. Write $\Ol_0 := \Zl_{(e)}[\theta_n]$ and let $\Ol$ be the integral closure of $\Ol_0$ in 
$\Ql(\theta_n)$. Let 
\[
\begin{array}{rcl}
\Ol_0^+ & := & \{ \xi\in \Ql(\theta_n)\; :\; \Tr_{\Ql(\theta_n)|\Ql}(\xi \Ol_0)\tm\Zl_{(e)}\} \\
\Ol^+   & := & \{ \xi\in \Ql(\theta_n)\; :\; \Tr_{\Ql(\theta_n)|\Ql}(\xi \Ol)\tm\Zl_{(e)}\}\; .\\
\end{array}
\]
Then $\Ol_0\tm\Ol\tm\Ol^+\tm\Ol_0^+$. But $\Ol_0 = \Ol_0^+$, since the $\Zl_{(e)}$-linear determinant of this embedding is given by the discriminant 
$\Delta_{\Zl[\theta_n]|\Zl}$, which is a unit in $\Ol_0$.
\end{Lemma}

We resume.

\begin{Proposition}[\cite{Ca38},\cite{Ha74}, cf.\ {\cite[\rm p.\ 115]{Go83}}]
\label{PropRes}
The extension $\Ql(\theta_n)|\Ql$ is galois of degree \linebreak[4] $[\Ql(\theta_n):\Ql] = (q-1)q^{n-1}$, with Galois group isomorphic to $(\Zl/f^n)^\ast$. The 
integral closure of $\Zl$ in $\Ql(\theta_n)$ is given by $\Zl[\theta_n]$. We have $\Zl[\theta_n]\theta_n^{[\Ql(\theta_n):\Ql]} = \Zl[\theta_n]f$. In particular, 
$\theta_n$ is a prime element of $\Zl[\theta_n]$, and the extension $\Zl_{(f)}[\theta_n]|\Zl_{(f)}$ of discrete valuation rings is purely ramified.
\end{Proposition}

\subsection{Coefficient valuation bounds}
\label{SubSecCVB}

Denote $\Fl_n = \Ql(\theta_n)$. Let $\El_n = \mb{Fix}_{C_{q-1}} \Fl_n$, so $[\El_n:\Ql] = q^{n-1}$. Let 
\[
\varpi_n \; =\; \Nrm_{\Fl_n|\El_n}(\theta_n) \; =\; \prod_{e\in (\Zl/f)^\ast} \theta_n^{e^{q^{n-1}}}\; .
\]
The minimal polynomial $\mu_{\theta_n,\Fl_{n-1}}(X) = P_f(X) - \theta_{n-1}$ together with $X\, |\, P_f(X)$ shows that 
$\Nrm_{\Fl_n|\Fl_{n-1}}(\theta_n) = \theta_{n-1}$, whence $\Nrm_{\El_n|\El_{n-1}}(\varpi_n) = \varpi_{n-1}$. Note that 
$\varpi_1 = \prod_{e\in (\Zl/f)^\ast} \theta_1^e = \Psi_{\! f}(0) = f$.

The extension $\Zl_{(f)}[\varpi_n]$ is a discrete valuation ring with maximal ideal generated by $\varpi_n$, purely ramified of index $q^{n-1}$ over $\Zl_{(f)}$. 
In particular, $\El_n = \Ql(\varpi_n)$. 

\begin{footnotesize}
\begin{quote}
\begin{Example}
\label{ExCar5}\rm
Let $r = 3$ and $f(Y) = Y^2 + 1$, so $q = 9$. A Magma calculation shows that
\[
\begin{array}{l}
\varpi_2 = \theta_2^{60} - Y\theta_2^{58} + Y^2\theta_2^{56} + (-Y^9 \! -\! Y^3 \! -\! Y)\theta_2^{42} + (Y^{10} \! +\! Y^4 \! +\! Y^2 \! +\! 1)\theta_2^{40} \\
+\; (-Y^{11} \! -\! Y^5 \! -\! Y^3 \! +\! Y)\theta_2^{38} + (-Y^6 \! -\! Y^4 \! -\! Y^2)\theta_2^{36} + (Y^7 \! +\! Y^5 \! +\! Y^3 \! +\! Y)\theta_2^{34} \\
+\; (-Y^8 \! -\! Y^6 \! +\! Y^4 \! -\! Y^2 \! -\! 1)\theta_2^{32} 
    + (-Y^5 \! +\! Y^3 \! -\! Y)\theta_2^{30} + (Y^{18} \! -\! Y^{12} \! -\! Y^{10} \! +\! Y^6\!\! -\! Y^4 \! +\! Y^2)\theta_2^{24} \\
+\; (-Y^{19} \! +\! Y^{13} \! +\! Y^{11} \! +\! Y^9 \! -\! Y^7 \! +\! Y^5 \! +\! Y)\theta_2^{22} \\
+\; (Y^{20} \! -\! Y^{14} \! -\! Y^{12} \! +\! Y^{10} \! +\! Y^8 \! -\! Y^6 \! -\! Y^4 \! +\! Y^2 \! +\! 1)\theta_2^{20} \\
+\; (-Y^{15} \! -\! Y^{13} \! -\! Y^{11} \! -\! Y^9 \! +\! Y^7 \! +\! Y^5 \! -\! Y^3)\theta_2^{18} 
    + (Y^{16} \! +\! Y^{14} \! +\! Y^{12} \! -\! Y^{10} \! -\! Y^8 \! -\! Y^2)\theta_2^{16} \\
+\; (-Y^{17} \! -\! Y^{15} \! +\! Y^{13}\! +\! Y^{11} \! +\! Y^7 \! +\! Y^5 \! -\! Y^3 \! +\! Y)\theta_2^{14} \\
+\; (-Y^{14} \! -\! Y^{12} \! +\! Y^{10} \! -\! Y^8 \! -\! Y^6 \! -\! Y^4 \! +\! Y^2 \! +\! 1)\theta_2^{12} 
    + (-Y^{13} \! +\! Y^{11} \! -\! Y^7 \! +\! Y^3)\theta_2^{10} \\
+\; (Y^{14} \! -\! Y^{12} \! -\! Y^{10} \! +\! Y^6 \! +\! Y^4)\theta_2^8 + (-Y^{11} \! -\! Y^7 \! +\! Y^5 \! +\! Y^3 \! +\! Y)\theta_2^6 + 
(Y^8 \! +\! Y^6 \! +\! Y^2 \! +\! 1)\theta_2^4\; . \\
\end{array}
\]
With regard to section \ref{CaseY}, we remark that $\varpi_2\neq \pm\,\theta_2^{q - 1}$.
\end{Example}
\end{quote}
\end{footnotesize}

\begin{Lemma}
\label{Lem10a}
For all $n\geq 2$, we have $\;\varpi_n^q\;\con_{\varpi_n^{q-1}f}\;\varpi_{n-1}\;$.

\rm
We claim that $\theta_n^q\con_{\theta_n f}\theta_{n-1}$. In fact, the non-leading coefficients of the Eisenstein polynomial $\Psi_{\! f}(X)$ are divisible by 
$f$, so that the congruence follows by $\theta_{n-1} - \theta_n^q = P_f(\theta_n) - \theta_n^q = \theta_n(\Psi_{\! f}(\theta_n) - \theta_n^{q - 1})$. Letting 
$\w T = \Zl_{(f)}[\theta_n]$ and $(\w t,\w s,t,s) = (\theta_n,\theta_{n-1},\varpi_n,\varpi_{n-1})$, (\ref{Lem9}) shows that $1 - \theta_n^q/\theta_{n-1}$ 
divides $1 - \varpi_n^q/\varpi_{n-1}$. Therefore, $\theta_n f\theta_{n-1}^{-1}\varpi_{n-1}\; |\;\varpi_{n-1} - \varpi_n^q$.
\end{Lemma}

Now suppose given $m\geq 1$. To apply (\ref{Prop7}), we let $R_i = \Zl_{(f)}[\varpi_{m+i}]$ and $r_i = \varpi_{m+i}$ for $i\geq 0$. We continue to denote
$$
\begin{array}{r}
\mu_{\varpi_{m+i},\,\El_m}(X) \; = \; \mu_{r_i,\,K_0}(X)\; =\; X^{q^i} + \Big(\sum_{j\in [1,q^i - 1]} a_{i,j} X^j\Big) - \varpi_m \vspace*{1mm}\\
\; \in\; R_0[X] \; =\; \Zl_{(f)}[\varpi_m][X]\; , \hspace*{16mm} \\
\end{array}
\leqno (\#)
$$
and $v_q(j) = \max\{ \alpha\in\Z_{\geq 0}\; :\; j\con_{q^\alpha} 0 \;\}$.

\begin{Theorem}
\label{Th11a}\Absit
\begin{itemize}
\item[{\rm (i)}] We have $f^{i - v_q(j)}\; |\; a_{i,j}$ for $i\geq 1$ and $j\in [1,q^i-1]$. 
\item[{\rm (i$'$)}] If $j < q^i (q-2)/(q-1)$, then $f^{i - v_q(j)} \varpi_m\; |\; a_{i,j}$.
\item[{\rm (ii)}] We have $a_{i,j} \con_{f^{i+1}} a_{i+\be,q^\be j}$ for $i\geq 1$, $j\in [1,q^i-1]$ and $\be\geq 1$. 
\item[{\rm (ii$'$)}] If $j < q^i (q-2)/(q-1)$, then $a_{i,j} \con_{f^{i+1}\varpi_m} a_{i+\be,q^\be j}$ for $\be\geq 1$.
\end{itemize}

\rm
Assumption (\ref{Ass6_5}) is fulfilled by virtue of (\ref{Lem10a}), whence the assertions follow by (\ref{Prop7}).
\end{Theorem}

\subsection{Some exact valuations}

Let $m\geq 1$ and $i\geq 0$. We denote $R_i = \Zl_{(f)}[\varpi_{m+i}]$, $r_i = \varpi_{m+i}$, $K_i = \fracfield R_i$, \linebreak 
$\w R_i = \Zl_{(f)}[\theta_{m+i}]$ and $\w r_i = \theta_{m+i}$. We obtain $\Dfk_{\w R_i|\w R_0} = (f^i)$ and 
$\Dfk_{\w R_i|R_i} = (\w r_i^{q-2})$ {\bf\cite[\rm III.\S3, prop.\ 13]{Se62}}, whence
$$
\Dfk_{R_i|R_0} \= \left(\mu'_{r_i,K_0}(r_i)\right) \= \left(f^i r_i^{q^i - 1 - (q^i - 1)/(q - 1)}\right)\; .
\leqno (\ast\ast)
$$
Therefore, $f^i r_i^{q^i - 1 - (q^i - 1)/(q - 1)}$ divides $j a_{i,j} r_i^{j-1}$ for $j\in [1,q^i - 1]$, which is an empty assertion if $j\con_p 0$.
Thus (\ref{Th11a}.\,i,\,i$'$) do not follow entirely.

However, since only for $j = q^i - (q^i - 1)/(q - 1)$ the valuations at $r_i$ of $f^i r_i^{q^i - 1 - (q^i - 1)/(q - 1)}$ and $j a_{i,j} r_i^{j-1}$
are congruent modulo $q^i$, we conclude by ($\ast\ast$) that they are equal, i.e.\ that
\[
a_{i,q^i - (q^i - 1)/(q - 1)}\,/\, f^i\;\in\; R_0^\ast\; .
\]

\begin{Corollary}
\label{CorDiff2}
We have
\[
a_{i,q^i - (q^i - q^\be)/(q - 1)}\,/\, f^{i-\be}\;\in\; R_0^\ast\= \Zl_{(f)}[\varpi_m]^\ast\Icm\mb{for $\;\;\be\in [0,i-1]$}\; .
\]

\rm
This follows by (\ref{Th11a}.ii) from what we have just said.
\end{Corollary}

\subsection{A simple case}
\label{CaseY}

Suppose that $f(Y) = Y$ and $m\geq 1$. Note that 
\[
\varpi_{m+1} \; =\; \prod_{e\in\sFu{q}^\ast} \theta_{m+1}^e \; =\; \prod_{e\in\sFu{q}^\ast} e\theta_{m+1} \; =\; -\theta_{m+1}^{q-1}\; . 
\]

\begin{Lemma}
\label{LemCar11}
We have
\[
\mu_{\varpi_{m+1},\El_m}(X)\; =\; - \varpi_m + \sum_{j\in [1,q]} Y^{q - j} X^j\; .
\]

\rm
Using the minimal polynomial $\mu_{\theta_{m+1},\Fl_m}(X) = P_Y(X) - \theta_m = X^q + YX - \theta_m$, we get
\[
\begin{array}{rl}
  & - \varpi_m + \sum_{j\in [1,q]} Y^{q - j} \varpi_{m+1}^j \\
= & \theta_m^{q - 1} + (Y^{q + 1} - \theta_{m+1}^{q^2 - 1})/(Y + \theta_{m+1}^{q - 1}) - Y^q \\
= & (Y\theta_m^{q - 1}\theta_{m+1} + \theta_m^{q - 1} \theta_{m+1}^q - \theta_{m+1}^{q^2} - Y^q \theta_{m+1}^q)/(\theta_{m+1}(Y + \theta_{m+1}^{q - 1}))  \\
= & 0\; . \\
\end{array}
\]
\end{Lemma}

\begin{Corollary}
\label{CorBU}
Let $m,i\geq 1$. We have
\[
\mu_{\varpi_{m+i},\El_m}(X) \;\con_{Y^2}\; X^{q^i} + YX^{(q-1)q^{i-1}} - \varpi_m\; .
\]

\rm
This follows from (\ref{LemCar11}) using (\ref{Th11a}.ii).
\end{Corollary}

\begin{quote}
\begin{footnotesize}
\begin{Remark}
\label{RemCar11_5}\rm
(\ref{LemCar11}) also holds if $p = 2$. 
\end{Remark}

\begin{Conjecture}
\label{ConjCar12}\rm
Let $m,i\geq 1$. We use the notation of $(\#)$ above, now in the case $f(Y) = Y$.
For $j\in [1,q^i]$, we write $q^i - j = \sum_{k \in [0,i-1]} d_k q^k$ with $d_k\in [0,q - 1]$. Consider the following conditions.
\begin{itemize}
\item[(i)] There exists $k\in [0,i-2]$ such that $d_{k+1} < d_k$. 
\item[(ii)] There exists $k\in [0,i-2]$ such that $v_p(d_{k+1}) > v_p(d_k)$. 
\end{itemize}
If (i) or (ii) holds, then $a_{i,j} = 0$. If neither (i) nor (ii) holds, then 
\[
v_{\varpi_m}(a_{i,j}) \; =\; q^{m-1}\cdot\sum_{k \in [0,i-1]} d_k\; .
\] 
\end{Conjecture}

\begin{Remark}
\label{RemCar13}\rm
We shall compare (\ref{CorDiff2}) with (\ref{ConjCar12}). If $j = q^i - (q^i - q^\be)/(q - 1)$ for some $\be\in [0,i-1]$, then 
$q^i - j = q^{i-1} + \cdots + q^\be$. Hence $\sum_{k\in [0,i-1]} d_k = i-\be$, and so according to (\ref{ConjCar12}), $v_{\varpi_m}(a_{i,j})$ should 
equal $q^{m-1}(i - \be)$, which is in fact confirmed by (\ref{CorDiff2}).
\end{Remark}
\end{footnotesize}
\end{quote}

%% file: wref.tex
\parskip0.0ex
\begin{footnotesize}

\parskip1.2ex

\vspace*{1cm}

\begin{flushright}
Matthias K\"unzer\\
Universit\"at Ulm\\
Abt.\ Reine Mathematik\\
D-89069 Ulm\\
kuenzer@mathematik.uni-ulm.de\\
\vspace*{1cm}
Eduard Wirsing\\
Universit\"at Ulm\\
Fak.\ f.\ Mathematik\\
D-89069 Ulm\\
ewirsing@mathematik.uni-ulm.de\\
\end{flushright}
\end{footnotesize}